\documentclass[12pt]{amsart}
\usepackage{amsmath}
\usepackage{amssymb}

\newcommand{\bi}{\bf\itshape}

\catcode`@=11 \@mparswitchfalse  

\numberwithin{equation}{section}   

\renewcommand{\(}{\left(}
\renewcommand{\)}{\right)}
\renewcommand{\[}{\left[}
\renewcommand{\]}{\right]}

\newcommand{\R}{{\mathbf R}}		
\newcommand{\e}{\varepsilon}            

\newcommand{\nn}{\nonumber}	   	
\newcommand{\note}[1]{}                    

\newcounter{mnotecount}[page]

\newcommand{\der}[2]{\frac{\partial #1}{\partial #2}}

\swapnumbers  

\newtheorem{thm}{Theorem}[section] 
\newtheorem{lemma}[thm]{Lemma}
\newtheorem{prop}[thm]{Proposition}
\newtheorem{cor}[thm]{Corollary}

\theoremstyle{definition}
\newtheorem{definition}[thm]{Definition}

\theoremstyle{remark}
\newtheorem{remark}[thm]{Remark}

\renewcommand{\:}{\colon}
\newcommand{\f}{\partial}

\renewcommand{\phi}{\varphi}
\newcommand{\image}{\operatorname{\rm image}\nolimits}
\newcommand{\ol}[1]{\overline{#1}}
\newcommand{\Span}{\operatorname{\rm span}\nolimits}
\newcommand{\bv}{{\mathbf v}}

\newcommand{\gt}{\theta}
\newcommand{\cd}{,\dots,}
\newcommand{\caxis}{C_{\text{\rm axes}}}
\newcommand{\eq}[1]{~(\ref{#1})}
\renewcommand{\setminus}{\smallsetminus}

\title[First Order Non-linear PDE in Two Variables]
{Solutions Near Singular Points to the Eikonal and Related First Order
Non-linear Partial Differential Equations in Two Independent
Variables}

\author[Cornea]{Emil Cornea${}^{*}$}
\address{Northern Illinois University
Mathematical Sciences
DeKalb, IL 60115-2888, USA}
\email{cornea\char'100math.niu.edu}
\thanks{${}^\dagger$Supported in part by ONR Grant No.
DAAH04-96-1-0326 and DoD Grant No. N00014-97-1-0806}
\author[Howard]{Ralph Howard${}^{\dagger}$}
\address{Department of Mathematics
University of South Carolina
Columbia, S.C. 29208, USA}
\email{howard\char'100math.sc.edu}
\thanks{${}^{\dagger}$Supported in part by DoD Grant No. N00014-97-1-0806}
\author[Martinsson]{Per-Gunnar Martinsson${}^{\ddagger}$}
\address{Texas Institute for Computational and Applied Mathematics, MC C0200,
University of Texas-Austin, Austin TX 78712, USA
}
\email{pgm\char'100ticam.utexas.edu}
\thanks{${}^{\ddagger}$Supported in part by ONR Grant Nos. N00014-97-11059 and
N00014-94-11163}

\begin{document}

\begin{abstract}  A detailed study of solutions to the
first order partial differential equation $H(x,y,z_x,z_y)=0$, 
with special emphasis on the eikonal equation $z_x^2+z_y^2=h(x,y)$, 
is made near points where the equation becomes singular in the 
sense that $dH=0$, in which case the method of characteristics does not apply.  
The main results are that there is a strong lack of uniqueness of solutions near
such points and that solutions can be less regular than both the
function $H$ and the initial data of the problem, but that this loss
of regularity only occurs along a pair of curves through the singular
point.  The main tools are symplectic geometry and the Sternberg
normal form for Hamiltonian vector fields.
\end{abstract}

\maketitle


\section{Introduction}
\label{sec:intro}

The eikonal equation in two independent variables
for $z = z(x,y)$ is
\begin{equation}\label{eqn:ikon}
z_x^2+z_y^2=h(x,y)
\end{equation}
where $h$ is a non-negative smooth function on the plane $\R^2$, and
subscripts denote partial derivatives.  
Near points $(x_0,y_0)$ where $h\neq0$ all local
solutions to this equation can be constructed by the method of
characteristics
(\emph{cf.}~\cite[Chap.~2]{Arnold:geoODE}, \cite[Chap.~10]{Spivak:comp5}) and
questions of local existence, uniqueness, and regularity are fully
understood.  However, near points $(x_0,y_0)$ where $h$ vanishes the
picture is much less complete.  Our main goal is to study these
questions for\eq{eqn:ikon} (and more generally for Hamilton-Jacobi
equations $H(x,y,z_x,z_y)=0$) in detail and show, contrary to some
assertions in the literature, that near a point where $h$ has zero as
a non-degenerate local minimum value there are
generally infinitely many local solutions.  The model case for such a
problem is
\begin{equation}\label{eq:toy}
z_x^2+z_y^2=a^2x^2+b^2y^2
\end{equation}
where $a$ and $b$ are positive.  For ease of statement we discuss
our general results in terms of this special case.  We normalize
the solution so that $z(0,0)=0$.  Then\eq{eq:toy} has at least the four
solutions $z_{\pm,\pm}=\frac12(\pm ax^2\pm by^2)$. If $a\ne b$ then
(see Proposition~\ref{lemma:sec} below) any $C^3$ solution near $(0,0)$
is of the form
$$
z(x,y)=\frac12(\pm ax^2\pm by^2)+O((|x|+|y|)^3).
$$
The two solutions $z_{+,+}=\frac12(ax^2+ by^2)$ and
$z_{-,-}=-\frac12(ax^2+ by^2)$ are unique in the sense that any $C^2$
solution $\tilde{z}$ with leading terms $\tilde{z}=\frac12(ax^2+
by^2)+O((|x|+|y|)^3)$ will agree with $z_{+,+}$ in a neighborhood of
$(0,0)$ (and on all of $\R^2$ if $\tilde{z}$ is globally defined),
with a similar statement holding for $z_{-,-}$.
(See~\cite{Bruss:eikonal} or Theorem~\ref{Bruss} below for the local
result and Oliensis~\cite{Oliensis:Uniqueness} 
for the global
version).  However, 
the other two solutions are far from unique.  In
Section~\ref{sec:basic} we show that for the solution
$z_{+,-}=\frac12(ax^2-by^2)$ there is an infinite dimensional family
of $C^\infty$ solutions of the form
$\tilde{z}=\frac12(ax^2-by^2)+O((|x|+|y|)^3)$ 
that do not agree with $z_{+,-}$
in any neighborhood of the origin (Theorem~\ref{saddle}).  (This gives
counterexamples to some claims in the literature,
\emph{cf.}~\cite[p.~1094]{Saxberg}, which would imply that $z_{\pm,\pm}$ are
the only solutions.)  It is also shown that for each $k\ge 2$ there
are  solutions that are $C^k$, but not $C^{k+1}$.  The method of
proof is to use a variation of the well known method of
characteristics for first order equations.  Data is given along the
two curves
$$
c_{+}(t)=(t,t),\quad c_{-}(t)=(t,-t)
$$
(two curves rather than just one because of the nature of the
singularity) and if the data is $C^k$ along these curves and agrees
with the data for $z_{+,-}$ to order $l\le k$ at the
origin then there is a surprising regularity phenomenon.  The solution
will be $C^{k+1}$ on $\R^2\setminus
\{xy=0\}$, as expected from the method of characteristics, but when
considered as a function on all of $\R^2$ the regularity is only $C^n$
with
$n:=\lceil\min\{\frac{(l+1)a}{a+b},
\frac{(l+1)b}{a+b}\}\rceil-1$ (where $\lceil\cdot \rceil$ is the
ceiling function) and
so there is a drop of regularity along the coordinate axes.  In
general, if the data is $C^k$ then the solution $\tilde{z}$ is $C^{k+1}$
off the coordinate axes, but on the coordinate axes the regularity is
determined by the order of contact of $\tilde{z}$ with the
``standard'' solution $z_{+,-}$.

This existence, lack of uniqueness, and jump of regularity along a
pair of distinguished curves is not  specific to the equation\eq{eq:toy}
$z_x^2+z_y^2=a^2x^2+b^2y^2$ but generalizes to equations
$z_x^2+z_y^2=h(x,y)$ where $h$ is of the form
\begin{equation}\label{gen:h}
h(x,y)=a^2x^2+b^2y^2+O((|x|+|y|)^3)
\end{equation}
(and more generally equations of the form $H(x,y,z_x,z_y)=0$) under
the extra assumption that the numbers $a$ and $b$ are linearly
independent over the rational numbers.  The proof is based on using
the Sternberg normal form for such an equation (and the existence of
this normal form is only guaranteed  when $a$ and $b$ are linearly
independent over the rationals) to reduce the calculations to
manageable proportions.

The condition that
 $a$ and $b$ are linearly independent over the rational
 numbers is exactly the condition to insure that\eq{gen:h} has
 a solution in formal power series with leading terms
 $z=\frac12(ax^2-by^2)+O((|x|+|y|)^3)$.  We include these calculations
 and use them to show that if $a$ and $b$ are linearly dependent over
 the rationals and $m$ and $n$ are integers such 
 that $ma-nb=0$ and
 $m+n\ge4,$ then $z_x^2+z_y^2=ax^2+by^2+x^my^n$ has no $C^k$ solution
 for $k\ge m+n$ with leading terms
 $z=\frac12(ax^2-by^2)+O((|x|+|y|)^3)$.  Thus the independence
 condition on $a$ and $b$ is both necessary and sufficient for the
 existence of smooth saddle type solutions.

A secondary goal of this paper is to 
advocate the use of
differential geometric methods, especially for differential forms,
symplectic geometry, and normal forms such as the Sternberg normal
form, in working with first order equations.  Thus we have included
some expository material in Section~\ref{sec:sym} about symplectic
geometry and its application to the method of characteristics for first
order equations.

The article
can be summarized as follows.  In Section~\ref{sec:sym} we present
the basic methods of symplectic geometry as applied to the method 
of characteristics for first order equations of the form
$H(x,y,z_x,z_y)=0$ in two independent variables.  For several reasons
the theory is easier in two dimensions and our hope is that this will
be useful in understanding the method in a concrete setting.  This
section also gives a proof of a slight generalization of a Theorem of
Bruss~\cite{Bruss:eikonal} on the existence of concave and convex
solutions to $z_x^2+z_y^2=h(x,y)$.

In Section~\ref{sec:main} we first give a detailed discussion of the
model equation $z_x^2+z_y^2=a^2x^2+b^2y^2$.  This leads to the
existence, non-uniqueness and lack of regularity results already
mentioned.  This is followed by a generalization of these results to
equations $H(x,y,z_x,z_y)=0$.  The same type of existence,
non-uniqueness and lack of regularity holds, but we were not able to
give quite as precise an analysis of the regularity.  The proofs make
essential use of both symplectic geometry and the Sternberg normal
form.

Section~\ref{sec:formal} looks at solutions to $z_x^2+z_y^2=h(x,y)$
from the point of view of formal power series.  This leads to examples
of equations that do not have any smooth solutions of saddle type.

In an appendix we include precise statements of two of the geometric
tools we use: the stable submanifold theorem and the Sternberg normal
form.

Finally we mention that equations of the type $z_x^2+z_y^2=h(x,y)$ and
more generally $H(x,y,z_x,z_y)=0$ have been of interest in computer
vision and related fields because of the ``shape from shading
problem'', where the goal is to reconstruct a surface $z(x,y)$
given a gray-scale image of it. If the surface is matte then the
intensity of the reflected light $I$ can be modeled as being
proportional to the scalar product of the surface normal and the
illumination direction, $I \sim \hat{n} \cdot \hat{v}$. If the viewing
and the illumination directions coincide, and we choose this direction
as the $z$-direction, then $\hat{n} = (-z_{x},\, -z_{y},\, 1)/\sqrt{1
+ z_{x}^{2} + z_{y}^{2}}$ and we get, after scaling,
$$
I(x,y) = \hat{n} \cdot \hat{z} = \frac{1}{\sqrt{1 + z_{x}^{2} + z_{y}^{2}}},
$$
and equation (\ref{eqn:ikon}) can be recovered by setting $h(x,y) =
1/I(x,y)^{2} - 1$.  Note that in applications the data function
typically has quite low regularity.  In this case using viscosity
solutions has proven to be effective, see \cite{lions:SFS} and
\cite{rouy:SFS}.


\section{Application of symplectic geometry to the eikonal equation}
\label{sec:sym}

\subsection{Review of the method of characteristics and construction
of concave and convex solutions.}  Let $\R^4$ have coordinates
$x,y,p,q$. Then the {\bi symplectic form\/} on $\R^4$ is
$$
\omega:=dp\wedge dx+dq\wedge dy.
$$
Let $N^2\subset \R^4$ be an imbedded surface.  Then $N^2$
{\bi projects\/} on an open set $U\subset\R^2$ iff there  
are continuous functions $p(x,y)$ and $q(x,y)$ so that
$$
N^2=\{(x,y,p(x,y),q(x,y)): (x,y)\in U\}.
$$
The submanifold $N^2$ is a {\bi jet of a function\/} iff there is an
open set $U\subseteq \R^2$ and a function $z \in C^{1}(U)$ such
that
$$
N^2=\{(x,y,z_x(x,y),z_y(x,y)): (x,y)\in U\}.
$$
It is clear that if $N^2$ is a jet of a function, then it projects,
but the converse is not true.  The following is a standard result, but
we include the short proof for those not familiar with the
differential geometric set up.

\begin{prop}  \label{prop:jet}
Assume that $N^2\subset \R^4$ is a simply connected two dimensional
submanifold of $\R^4$ of smoothness class $C^k$ for some $k\ge 1$.
Then $N^2\subset\R^4$ is the jet of some function (which will be
unique up to an additive constant) if and only if it projects over
some open set $U\subseteq\R^2$ and the restriction of the symplectic
form to $N^2$ vanishes.  (That is, $\omega(X,Y)=0$ for all vectors
tangent to $N^2$.)  Moreover, if $N^2$ is of class $C^k$ and is the jet
of a function $z$, then $z$ is of class $C^{k+1}$.  
\end{prop}

\begin{proof}
Assume that $N^2$ projects over $U$ so that
$N^2=\{(x,y,p(x,y),q(x,y)): (x,y)\in U\}$.  As $N^2$ is simply
connected the same is true of $U$.  Using $x,y$ as coordinates on
$N^2$, we see that the restriction of the symplectic form to $N^2$ is
\begin{align*}
\omega&=dp\wedge dx +dq\wedge dy=(p_x\,dx+p_y\,dy)\wedge dx
	+(q_x\,dx+q_y\,dy)\wedge dy\\
	&=(q_x-p_y)\,dx\wedge dy.
\end{align*}
Therefore the restriction of $\omega$ to $N^2$ vanishes iff
$p_y=q_x$. But as $U$ is simply connected this is exactly the
condition that there is a function $z$ (unique up to an additive
constant) so that $p=z_x$ and $q=z_y$.  If $N^2$ is of class $C^k$
then $p(x,y)$ and $q(x,y)$ are $C^k$ functions of $(x,y)$.  Thus $z$
is of class $C^{k+1}$.~\end{proof}

A two dimensional surface $N^2\subset \R^4$ is a {\bi Lagrangian surface\/}
iff the restriction of $\omega$ to $N^2$ vanishes.

Given a function $H\:\R^4\to \R$ the {\bi characteristic vector
field\/} of $H$ (also called the {\bi symplectic gradient\/}) is the
unique vector yield $\xi_H$ on $\R^4$ such that for all vectors $X$
$$
\omega(\xi_H,X)=-dH(X).
$$
Then a simple calculation yields that
\begin{equation}\label{char:form}
\xi_H=\der{H}{p}\der{}{x}+\der{H}{q}\der{}{y}
	-\der{H}{x}\der{}{p}-\der{H}{y}\der{}{q}.
\end{equation}
An immediate consequence of the definition of $\xi_H$ is that
$$
dH(\xi_H)=-\omega(\xi_H,\xi_H)=0
$$
and so $H$ is constant on the integral curves of $\xi_H$.

The following is the differential geometric justification of the
method of characteristics for the Hamilton-Jacobi equation
$H(x,y,z_x,z_y)=0$. Given a smooth function $H\:\R^4\to \R$ and $P\in
\R^4$,
 the restriction of $dH$ to $T(\R^4)_P$ is denoted $dH_P$.

\begin{prop}\label{soln}
Let $N^2$ be a simply connected two dimensional submanifold of $\R^4$.
We also assume that the following non-degeneracy condition holds:
\begin{enumerate}
\item[(ND)] The set of points $P\in N^2$ where $dH_P\ne 0$ is dense in $N^2$
(here $dH_P$ is being viewed as a linear functional on $\R^4$ and not
just restricted to $T(N^2)$).
\end{enumerate}
Then $N^2$ is the jet of a solution to the equation
$H(x,y,z_x,z_y)=0$ if and only if the
following three conditions hold
\begin{enumerate}
\item $N^2\subset \{P: H(P)=0\}$, 
\item $N^2$ projects over some open set $U\subseteq \R^2$, and
\item the characteristic vector field $\xi_H$ is tangent to $N^2$ at
every point of $N^2$.
\end{enumerate}
If $N^2$ satisfies these conditions and is a $C^k$ submanifold, then
the solution $z$ is of class $C^{k+1}$.
\end{prop}

\begin{remark}  For our future applications it is important to realize
that although $H$ is a smooth function, the zero set $\{H=0\}$ need
not be a smooth submanifold of $\R^4$.  For $a,b>0$ we will be
interested in $H(x,y,p,q)=p^2+q^2-a^2x^2-b^2y^2$.  The zero set
$\{H=0\}$ is then a cone that is singular at $(0,0,0,0)$.  However, if
$N^2$ is a smooth surface in $\{H=0\}$ that is everywhere tangent to
$\xi_H$ and which projects onto an open set in $\R^2$, then $N^2$ is
the jet of a solution to $z_x^2+z_y^2=a^2x^2+b^2y^2$.  This is because
$dH$ only vanishes at one point and thus the set of points on $N^2$
where $dH\ne 0$ is dense.\qed
\end{remark}

It is useful to give a name to submanifolds that satisfy two of the
conditions of the proposition:

\begin{definition}\label{invar-Lagrange}
Let $H \: \R^4\to \R$ be a $C^k$ function with $k\ge 2$ and
characteristic vector field $\xi_H$. Then a connected two dimensional
submanifold $N^2$ of class $C^1$ is an
{\bi invariant Lagrangian surface\/} iff
\begin{enumerate}
\item $N^2$ is a Lagrangian surface (\emph{i.e.} the restriction of $\omega$
to $N^2$ vanishes), and 
\item The characteristic vector field $\xi_H$ is tangent to $N^2$ at
all points of $N^2$.\qed
\end{enumerate}
\end{definition}

\begin{proof}[Proof of Proposition~\ref{soln}]
First assume that the three conditions hold.  Then in light of
Proposition~\ref{prop:jet} it is enough to show that the restriction
of $\omega$ to $N^2$ is zero.  If we are at a point $P\in N^2$ where
$dH\ne 0$, then from the definition of $\xi_H$ or the
formula~(\ref{char:form}) it follows that $\xi_H\ne0$.  Let $X$ be a
tangent vector to $N^2$ at $P$ linearly independent from $\xi_H(P)$.
Then $\{\xi_H(P),X\}$ is a basis for the tangent space
$T(N^2)_P$. Using the definition of $\xi_H$ and that $dH(X)=0$ (as
$H\big|_{N^2}=0$)
$$
\omega(\xi_H,X)=-dH(X)=0.
$$
Thus $\omega$ restricted to $N^2$ vanishes at $P$.  But we are
assuming that the set of points $P$ where $dH$ does not vanish is dense,
thus by continuity the restriction of $\omega$ to $N^2$ is zero on all
of $N^2$.\smallskip

Conversely, if $N^2$ is the jet of a solution, then clearly
$N^2\subset \{P: H(P)=0\}$ and $N^2$ projects onto an open subset
of $\R^2$ and also $H\big|_{N^2}=0$.  Let $P\in N^2$.  Then
$\omega(X,Y)=0$ for all $X,Y\in T(N^2)_P$.  But a calculation shows
that if $Z$ is any vector tangent to $\R^4$ at $P$ with
$\omega(Z,X)=0$ for all $X\in T(N^2)_P$, then $Z\in T(N^2)_P$.  But for
$X\in T(N^2)_P$ we have $dH(X)=0$ so (as above)
$\omega(\xi_H,X)=-dH(X)=0$.  Thus $\xi_H\in T(N^2)_P$.  This completes
the proof.~\end{proof}

This result makes the geometry of the method of characteristics
clearer than the classical presentations.  Let $H$ be of class $C^k$
for $k\ge2$.  Then the formula~(\ref{char:form}) makes it clear that
the vector field $\xi_H$ is of class $C^{k-1}$.  Let $\Phi^H_t$ be the
flow of $\xi_H$.  That is, $\Phi_0^H(P)=P$ and $c(t):=\Phi^H_t(P)$ is
an integral curve of the vector field $\xi_H$.  Then (see
\cite[Thm~1,~p.~80]{Lang:manifolds} 
or~\cite[p.~230]{Arnold:ODE}) the
map $(t,P)\mapsto \Phi^H_t(P)$ is $C^{k-1}$.  Now let $c\:(a,b)\to
\R^4$ be a curve of class $C^l$ ($l\ge 1$) and let $\pi\:\R^4\to\R^2$
be the projection $\pi(x,y,p,q)=(x,y)$.  Assume
\begin{enumerate}
\item $H(c(s))\equiv0$,
\item $c$ projects over an imbedded\footnote{A curve $c$ is imbedded
iff it is the image of a smooth injective map $\gamma: I \to \R^2$,
with $I$ an interval in $\R$,  so
that the velocity vector $\gamma'$ never vanishes and so that the
topology of $\gamma[I]$ as a subset of $\R^{2}$ 
is the same as the topology induced by $\gamma$ 
(\emph{ i.e.} $\gamma$ is a homeomorphism).}
curve of $\R^2$ (that is $s\mapsto
\pi(c(s))$ is an imbedded curve in $\R^2$), and
\item at all points $c(s)$ the vectors
$\pi_*c'(s)$ and $\pi_*\xi_H(c(s))$ are linearly independent.  
\end{enumerate}
Then let $F(s,t):=\Phi^H_t(c(s))$.  By the implicit function theorem
for $r$ small enough the submanifold $N^2:=\{F(s,t):(s,t)\in
(a,b)\times (-r,r)\}$ will project over an open subset $U$ of $\R^2$.
Also, as $H$ is constant along the flow of $\xi_H$ and $H(c(s))=0$, we
have that $H(F(s,t))=H(\Phi^H_t(c(s)))=0$.  Thus $H\big|_{N^2}=0$.  By
construction $\xi_H$ is tangent to $N^2$ and so $N^2$ is the jet of a
solution $z$ to $H(x,y,z_x,z_y)=0$ by Proposition~\ref{soln}.  The
regularity of $N^2$ is $C^{\min\{k-1,l\}}$ and therefore the
regularity of $z$ is $C^{\min\{k,l+1\}}$.

\subsection{Solutions near critical points of $H$.}

We now look at the more interesting case of finding solutions near a
critical point, $P_0$, of $H$.  Assume that $dH_{P_0}=0$ and by adding a
constant to $H$ we can assume that $H(P_0)=0$.  Then near $P_0$ the
set $\{H=0\}$ need not be a submanifold of $\R^4$.  Assume that $P_0$
is a non-degenerate critical point so that the Hessian of $H$ at $P_0$
is non-singular.  To make the notation easier we assume that
$P_0=(0,0,0,0)$.  Then using the form~(\ref{char:form}) we see that
the linearization of the characteristic system for this vector field
at the origin is
\begin{equation}\label{L-def}
\frac{d}{dt}\[\begin{matrix}x(t)\\ y(t)\\ p(t)\\ q(t)\end{matrix}\]
=\[\begin{matrix}H_{xp}&H_{yp}&H_{pp}&H_{pq}\\
		H_{xq}&H_{yq}&H_{pq}&H_{qq}\\
		-H_{xx}&-H_{xy}&-H_{xp}&-H_{xq}\\
		-H_{xy}&-H_{yy}&-H_{yp}&-H_{yq}\end{matrix}\]
		\[\begin{matrix}x(t)\\ y(t)\\ p(t)\\q(t)\end{matrix}\]
	=L\[\begin{matrix}x(t)\\ y(t)\\ p(t)\\q(t)\end{matrix}\],
\end{equation}
where all the second partial derivative are evaluated at $(0,0,0,0)$
and this equation defines $L$. Let $\det(D^2H)$ be the determinant of
the Hessian at $(0,0,0,0)$ and let
$$
c_2:=2{{H}_{{ xy}}}\,{{H}_{{ pq}}} - 2{{H}_{
{ xq}}}\,{{H}_{{ yp}}} + {{H}_{{ yy}}}{{H}_{{ qq}}}
 - {H}_{{ yq}}^{2} + {{H}_{{ xx}}}{{H}_{{ pp}}} - {
H}_{{ xp}}^{2}.
$$
The characteristic polynomial of $L$ is then, using that
$\det(L)=\det(D^2H)$, (which can be seen by noting that $L=J\,D^2H$, where
$J=\text{\tiny $\[\begin{matrix}0&I_2\\-I_2&0\end{matrix}\]$}$,
so that $\det(L)=\det(J)\det(D^2H)=\det(D^2H)$),
$$
\det(\lambda I-L)=\lambda^4+c_2\lambda^2+\det(D^2H).
$$
Therefore the eigenvalues are
\begin{equation}\label{eig-val}
\pm\sqrt{\frac{-c_2+\sqrt{c_2^2-4\det(D^2H)}}{2}},\quad
\pm\sqrt{\frac{-c_2-\sqrt{c_2^2-4\det(D^2H)}}{2}}.
\end{equation}
Assuming that there is no eigenvalue with zero real part we see that
there are exactly two eigenvaules with positive real part and two with
negative real part. Then let
$N^2_+$ be the local stable manifold for the critical point at
$(0,0,0,0)$ and $N^2_-$ the local unstable submanifold.  (Loosely
$N^2_{\pm}$ is the set of points $P\in\R^4$
so that $\lim_{t\to\pm\infty}\Phi^H_t(P)=(0,0,0,0)$.  See
Appendix~\ref{app:stable}.) Because of the condition on the
eigenvalues of $L$ both $N^2_+$ and $N^2_-$ are two dimensional.  If
$H$ is of class $C^k$ for $k\ge2$ then the flow $\Phi^H_t$ is of class
$C^{k-1}$.  Therefore $N^2_+$ and $N^2_-$ are $C^{k-1}$ submanifolds
of $\R^4$ (see \cite[Thm~5.20, p.~49]{Shub:stable}).  Then for a point
$P\in N^2_{\pm}$ we have, by the invariance of $H$ under the flow,
that $H(P)=\lim_{t\to\pm\infty}H(\Phi^H_t(P))=H(0,0,0,0)=0$.  Thus
$H\big|_{N^2_{\pm}}=0$.  Also it is clear that the characteristic
vector field $\xi_H$ is tangent to $N^2_{\pm}$.  The condition of
Proposition~\ref{soln} that does not hold automatically for
$N^2_{\pm}$ is that of being locally projectable.  Summarizing:

\begin{thm}\label{thm:stable}  Let $H\:\R^4\to\R$ be a $C^k$ function
with $k\ge2$.  Assume that $H$ has a non-degenerate critical point at
$(0,0,0,0)$ with $H(0,0,0,0)=0$ and assume that no eigenvalue of $L$
as defined above has zero real part.  Then there are exactly two
eigenvalues $\lambda_1,\lambda_2$ of $L$ that have positive real part.
Let $e_1$ and $e_2$ be the eigenvectors corresponding to $\lambda_1$
and $\lambda_2$, let ${\mathcal N}$ be the subspace spanned by $e_1$
and $e_2$ and let $\pi\:\R^4\to\R^2$ be the projection
$\pi(x,y,p,q)=(x,y)$.  Assume that the restriction $\pi\big|_{\mathcal
N}\:{\mathcal N}\to \R^2$ is nonsingular.  Then near $(0,0,0,0)$ the
unstable submanifold $N_-^2$ of $\xi_H$ is the jet of a $C^{k}$
solution $z$ to $H(x,y,z_x,z_y)=0$. An analogous statement holds
for the stable submanifold $N^2_+$ of $\xi_H$.
\end{thm}

\begin{proof}  The tangent space to $N^2_-$ at $(0,0,0,0)$
is ${\mathcal N}$ so if $\pi\big|_{\mathcal N}\:{\mathcal N}\to \R^2$
is nonsingular the implicit function theorem implies that near
$(0,0,0,0)$ the submanifold $N^2_-$ will project onto an open subset
of $\R^2$.  As noted all the other hypothesis of Proposition~\ref{soln}
hold.  This completes the proof.~\end{proof}

As an example of this assume that $H(x,y,p,q)=f(p,q)-h(x,y)$ where $f$
and $h$ have critical points at $(0,0)$ and the Hessians of $f$ and
$h$ at $(0,0)$ are both positive definite.  (This is a case that comes
up in the shape from shading problem.) Then the linear map $L$ is
$$
L=\[\begin{matrix}0&0&f_{pp}&f_{pq}\\ 0&0&f_{pq}&f_{qq}\\ h_{xx}&h_{xy}&0&0\\
	h_{xy}&h_{yy}&0&0\end{matrix}\]=
\[\begin{matrix}0&A\\ B&0\end{matrix}\],
$$
where $A$ and $B$ are positive definite  $2\times 2$ matrices.
As $A$ is positive definite it has a square root $A^{\frac12}$.
The matrix $A^{\frac12}BA^{\frac12}$ will also be positive definite
and  will consequently have positive eigenvalues.  Let $a^2$ and
$b^2$ be the eigenvalues of $A^{\frac12}BA^{\frac12}$, where
$a,b>0$, and ${\mathbf e}_{1}, {\mathbf e}_{2}$ the corresponding eigenvalues,
$$
A^{\frac12}BA^{\frac12}{\mathbf e}_1=a^2 {\mathbf e}_1,\quad
	A^{\frac12}BA^{\frac12}{\mathbf e}_2=b^2{\mathbf e}_2,
$$
and define
$$\begin{array}{ll}
{\mathbf v}_1:=\[\begin{matrix}A^{\frac12}{\mathbf e}_1\\ a
	A^{-\frac12}{\mathbf e}_1\end{matrix}\],&\quad 
{\mathbf v}_2:=\[\begin{matrix}A^{\frac12}{\mathbf e}_2\\ b
	A^{-\frac12}{\mathbf e}_2\end{matrix}\],\\ &\\
{\mathbf v}_3:=\[\begin{matrix}A^{\frac12}{\mathbf e}_1\\ -a
	A^{-\frac12}{\mathbf e}_1\end{matrix}\],&\quad 
{\mathbf v}_4:=\[\begin{matrix}A^{\frac12}{\mathbf e}_2\\ 
-b A^{-\frac12}{\mathbf e}_2\end{matrix}\].\end{array}
$$
By direct calculation, we find
$$
L{\mathbf v}_1=a {\mathbf v}_1,\quad L{\mathbf v}_2=b {\mathbf v}_2,\quad
L {\mathbf v}_3=-a {\mathbf v}_3,\quad L{\mathbf v}_4=-b {\mathbf v}_4.
$$
Thus the eigenvalues of $L$ are $\pm a$ and $\pm b$.  The eigenvectors
corresponding to the two eigenvalues with positive real part are
${\mathbf v}_1$ and ${\mathbf v}_2$.  For future reference we record some
elmentary facts about these eigenvalues and eigenvectors.  The proof
is left to the reader.

\begin{lemma} \label{lem:facts}
If $a\ne b$ then there are exactly six two dimensional subspaces of
the tangent space $T(\R^4)$ that are invariant under $L$,
corresponding to the six ways of choosing two of the eigenvectors
${\mathbf v}_1,\dots, {\mathbf v}_4$.  As $\omega({\mathbf
v}_1,{\mathbf v}_3)\ne 0$ and $\omega({\mathbf v}_2,{\mathbf v}_4)\ne
0$ the subspaces $\Span({\mathbf v}_1,{\mathbf v}_3)$ and
$\Span({\mathbf v}_2,{\mathbf v}_4)$ can never be tangent to jets of
functions and therefore can be disregarded from our considerations.
Each of the remaining four pairs $\{{\mathbf v}_1,{\mathbf v}_2\}$,
$\{{\mathbf v}_3,{\mathbf v}_4\}$, $\{{\mathbf v}_1,{\mathbf v}_4\}$,
and $\{{\mathbf v}_2,{\mathbf v}_3\}$ spans
a two dimensional subspace of the tangent
space $T(\R^4)_0$ on which $\omega$ vanishes.  If $\{{\mathbf v}_i,
{\mathbf v}_j\}$ is any one of these four pairs and $\pi\:\R^4\to
\R^2$ is the projection $\pi(x,y,p,q)=(x,y)$ then $\pi_*{\mathbf v}_i$
and $\pi_*{\mathbf v}_j$ are
linearly independent.  Therefore, if $N^2$ is any two dimensional
submanifold tangent to $\Span({\mathbf v}_i,{\mathbf v}_j)$ then, by
the implicit function theorem, $N^2$ locally projects over some open
neighborhood $U$ of $(0,0)$ in $\R^2$. Finally, $\Span({\mathbf
v}_1,{\mathbf v}_2)$ is the tangent space at the origin to the
unstable submanifold of $\xi_H$ and $\Span({\mathbf v}_3,{\mathbf
v}_4)$ is the tangent space at the origin to the stable submanifold
of~$\xi_H$.~\qed
\end{lemma}

The following gives a generalization of a theorem of
Bruss~\cite{Bruss:eikonal}, who considered the special case of
$u_x^2+u_y^2=h(x,y)$.  The proof uses the same geometric idea as
Bruss's proof.

\begin{thm}[A.~Bruss~\cite{Bruss:eikonal}]\label{Bruss}
Let $f(p,q)$ and $h(x,y)$ be functions of class $C^k$, with $k\ge2$,
and assume that $f$ and $h$ both have critical points with positive
definite Hessians at $(0,0)$.  Then in the class of functions that are
concave near $(0,0)$ there is a unique (up to an additive constant)
solution to $f(z_x,z_y)=h(x,y)$.  This solution is $C^{k}$ and has as
its jet near $(0,0)$ the stable manifold of the characteristic vector
field $\xi_H$ of $H(x,y,p,q)=f(p,q)-h(x,y)$.  (Likewise, there is a
unique convex solution, it is $C^{k}$ and has as jet near $(0,0)$ the
unstable submanifold of $\xi_H$.)
\end{thm}

\begin{proof} Follows from the discussion above.
\end{proof}

\section{Construction of saddle type solutions near a regular critical
point}
\label{sec:main}

\subsection{Construction of all solutions to $z_x^2+z_y^2=a^2x^2+b^2y^2$.}
\label{sec:basic}

We will now construct all ``saddle point'' solutions of the equation
$$
z_x^2+z_y^2=a^2x^2+b^2y^2,
$$
where $a,b>0$.  The analysis here is a model for the more general
setup covered in Section~\ref{sec:general-case}.  This special case is
also of interest as it is possible to be somewhat more precise about
the regularity of the solutions.  What makes this equation especially
easy to analyze is that the components of the characteristic vector
field are linear so that finding the flow involves no more than linear
algebra.  We will assume that we have a solution defined near $(0,0)$
with $z(0,0)=0$.  If $a\ne b$, then by a formal power series argument
as in the proof of Proposition~\ref{lemma:sec} below, we have that if
$z \in C^3$ then near $(0,0)$ the first few terms of its Taylor series
are $z=\frac12(\pm ax^2\pm by^2)+O((|x|+|y|)^3)$ and if $a=b$ we can
bring $z$ to this form by a rotation of the axes.  If
$z=\frac12(ax^2+ by^2)+O((|x|+|y|)^3)$, then Theorem~\ref{Bruss}
applies, since the function is convex.  Hence $z=\frac12(ax^2+by^2)$ is
the unique solution.  Likewise, if $z=-\frac12(ax^2+
by^2)+O((|x|+|y|)^3)$, then $z=-\frac12(ax^2+ by^2)$.

Next we
look for a solution of the form
$$
z=\frac12(ax^2-by^2)+O((|x|+|y|)^3).
$$
We will use the Hamiltonian $H=\frac12(p^2+q^2-a^2x^2-b^2y^2)$.  The
characteristic vector field is then
$\xi_H=p\der{}{x}+q\der{}{y}+a^2x\der{}{p}+b^2y\der{}{q}$.  The integral
curves of this vector field satisfy the differential equation,
\emph{ cf.} (\ref{L-def}),
$$ \frac{d}{dt} 
\[ \begin{matrix} x \\ y \\ p \\ q \end{matrix} \] =
\[ \begin{matrix} 0 & 0 & 1 & 0 \\ 0 & 0 & 0 & 1 \\
a^{2} & 0 & 0 & 0 \\ 0 & b^{2} & 0 & 0 \end{matrix} \]
\[ \begin{matrix} x \\ y \\ p \\ q \end{matrix} \] = 
L \[ \begin{matrix} x \\ y \\ p \\ q \end{matrix} \].$$
The eigenvectors of the matrix $L$ are
\begin{equation}
\label{e-vectors}
{{\mathbf v}}_{1}=\[\begin{matrix} 1\\ 0\\ a\\ 0\end{matrix}\],\quad
{{\mathbf v}}_{2}=\[\begin{matrix} 0\\ 1\\ 0\\ -b\end{matrix}\],\quad
{{\mathbf v}}_{3}=\[\begin{matrix} 1\\ 0\\ -a\\ 0\end{matrix}\],\quad
{{\mathbf v}}_{4}=\[\begin{matrix} 0\\ 1\\ 0\\ b \end{matrix}\],
\end{equation}
which satisfy
$$
L{{\mathbf v}}_{1}=a{{\mathbf v}}_{1},\quad L{{\mathbf v}}_2=-b
{{\mathbf v}}_2,\quad L{{\mathbf v}}_3=
-a{{\mathbf v}}_3,\quad L{{\mathbf v}}_4=b{{\mathbf v}}_4.
$$
The jet of the function $z=\frac12(ax^2-by^2)$ is the span of
${\mathbf v}_1$ and ${\mathbf v}_2$. 

To construct an invariant Lagrangian surface through $(0,0,0,0)$ we
start with two curves $\gamma_\pm\:\R\to \R^4$ and generate the
surface by moving these curves by the flow of the characteristic
vector field. Because of the nature of the singularity, 
two curves rather than just one are required. Toward this end
let $\phi_+,\phi_-
\:\R\to \R$ be two functions so that
\begin{enumerate}
\item $\phi_\pm$ are $C^k$ functions for some $k\ge1$ and
\item $\phi_\pm$ vanish to order $l$ at the origin for some $l\le k$
(where $l$ need not be an integer).  Specifically, this means that
there are functions $\ol{\phi}_\pm\:\R\to \R$ so that
\begin{equation}\label{phi-bar}
\phi_\pm(s)=\ol{\phi}_\pm(s)s^l
\end{equation}
where $\ol{\phi}_\pm(s)$ are $C^k$ functions on $\R\setminus\{0\}$
and the derivatives $\frac{d^j}{ds^j}\ol{\phi}_\pm(s)$ are bounded on
$[-1,1]\setminus \{0\}$ for $0\le j\le l$.
\end{enumerate}

Define two curves $\gamma_+, \gamma_-:\R\to\R^4$ by
$$
\gamma_{\pm}(s)=s{\mathbf v}_1\pm s{\mathbf v}_2
	+\phi_\pm(s){\mathbf v}_3\mp\frac{a^2}{b^2}\phi_\pm(s){\mathbf v}_4. 
$$
Then
$$
H(\gamma_\pm(s))=0
$$
for all $s$.  Let $e^{tL}$ be the exponential of $L$ so that
for $P\in\R^4$
the integral curve of the characteristic vector field through $P$ is
$t\mapsto e^{tL}P$.  Define $F_\pm\:\R^2\to \R^4$ by
\begin{equation}\label{F-def}
F_\pm(s,t):=e^{tL}\gamma_\pm(s)=se^{at}{\mathbf v}_1\pm se^{-bt}
 	{\mathbf v}_2+
\phi_{\pm}(s)e^{-at}{\mathbf v}_3\mp\frac{a^2}{b^2}\phi_{\pm}(s)e^{bt}{\mathbf v}_4.
\end{equation}
We now consider $F_+(s,t)$ for $s>0$.  Do the change of variables
$(s,t)\mapsto (u,v)$ with $u,v>0$ given by
$$ \left \{ \begin{array}{l} u = s e^{a t} \\
                   v = s e^{-b t} \end{array} \right.
\quad
\left \{ \begin{array}{l} t = \dfrac{1}{a + b}\ln \dfrac{u}{v} \\
s = u^{b / (a+b)} v^{a / (a+b)} 
\end{array} \right. . 
$$
In these coordinates we have that the image of $F_+(s,t)$ with $s>0$
can be parameterized by
{\def\s{u^{\frac{b}{a+b}}v^{\frac{a}{a+b}}}
\begin{align*}
G(u,v)=&u{\mathbf v}_1+v{\mathbf v}_2\\
&\quad+u^{\frac{-a}{a+b}}v^{\frac{a}{a+b}}\phi_+(\s){\mathbf v}_3
-\frac{a^2}{b^2}u^{\frac{b}{a+b}}v^{\frac{-b}{a+b}}\phi_+(\s){\mathbf v}_4.
\end{align*}}
\begin{lemma}\label{reg}
Let $\phi$ be a function on $[0,\infty)$ of class $C^k$ that vanishes
to order $l\le k$ at $t=0$ and let $1>\alpha,\beta>0$. Then the function
$$
E(u,v):=\begin{cases}u^{\alpha-1}v^{\beta}\phi(u^\alpha v^\beta),&u,v>0\\ 
	0,&\text{otherwise}\end{cases}
$$
is $C^n$ on $\R^2$ for all $n< \min\{(l+1)\alpha-1, (l+1)\beta\}$.
Conversely, if $E(u,v)$ is of class $C^n$, then $\phi$ vanishes of
order $l$ at $t=0$ for all $l>\max\{(n+1)/\alpha -1,
n/\beta-1\}$.
\end{lemma}

\begin{proof}
The hypothesis on $\phi$ implies that $\phi(t)=t^l\psi(t)$ where
$\psi$ has continuous derivatives up to order $k$ on the open interval
$(0,\infty)$ and is bounded on the closed interval $[0,\infty)$.
Thus $E$ can be rewritten as
$E(u,v)=u^{(l+1)\alpha-1}v^{(l+1)\beta}\psi(u^\alpha v^\beta)$ for
$u,v>0$.  The result now follows by direct calculation of the
derivatives.~\end{proof}

This directly implies:

\begin{thm} \label{saddle}
Let $\phi_\pm\:[0,\infty)\to \R$ be $C^k$ functions that vanish to
order $l\le k$ at $0$.  Let $F_\pm$ be defined by
equation~(\ref{F-def}) and let $N^2(\phi_+,\phi_-)$ be the closure of
the union of the images of $F_+$ and $F_-$, that is,
$$
N^2(\phi_+,\phi_-):=\overline{\image(F_+)\cup \image(F_-)}.
$$  
Set $n:=\lceil\min\{\frac{(l+1)a}{a+b},
\frac{(l+1)b}{a+b}\}\rceil-1$,
where $\lceil\cdot \rceil$ is the ceiling function.  Then
$N^2(\phi_+,\phi_-)$ is a $C^{n-1}$ submanifold of $\R^4$ that is the
jet of a $C^{n}$ solution $z(x,y)$ to the equation
$z_x^2+z_y^2=a^2x^2+b^2y^2$.  Conversely, any $C^n$ solution of this
equation is of this form for unique $\phi_+$
and $\phi_-$.  Regardless
of the value of $l$, the function $z$ will be of class $C^{k+1}$ on
$\R^2\setminus \{xy=0\}$.  Thus there is a decrease in regularity of
solutions from $C^{k+1}$ to $C^n$ along the coordinate axes.
\end{thm}

\begin{proof} 
All but the statements about the decrease in regularity along the
coordinate axes follow from Lemma~\ref{reg}.  The two curves where the
regularity of the surface $N^2(\phi_+,\phi_-)$ drops are along the
lines $u\mapsto G(u,0)$ (which projects down onto the $x$-axis) and
$v\mapsto G(0,v)$ (which projects down onto the $y$-axis).  All other
points of $N^2(\phi_+,\phi_-)$ are $C^k$.  As in
Proposition~\ref{soln}, this implies that $z$ is $C^{k+1}$ off of the
coordinate axes.
\end{proof}

\subsection{Solutions to $H(x,y,z_x,z_y)=0$ near a critical point of $H$.}
\label{sec:general-case}

We now consider an equation $H(x,y,z_x,z_y)=0$ near a critical point
of $H$.  We assume that $H \in C^\infty$, $H(0,0,0,0)=0$ and $dH=0$ at
the origin.  Let $L$ be the linearization at the origin of the
characteristic system defined by~(\ref{L-def}) and assume that the
eigenvalues of $L$ are $a,b,-a,-b$ where $a,b>0$.  (By (\ref{eig-val})
the eigenvalues are of this form if they are real and nonzero.)  Let
$e_1, e_2,e_3, e_4$ be the eigenvectors for $a, b,-a,-b$ respectively.
We make the two assumptions
\begin{equation}\label{a-b-ind}
\text{$a$ and $b$ are linearly independent over the rational numbers,}
\end{equation}
and 
\begin{equation}\label{ee-proj}
\left\{\begin{array}{l}
\{\pi_*e_1,\pi_*e_2\},\ 
\{\pi_*e_3,\pi_*e_4\},\  \{\pi_*e_1,\pi_*e_4\},\  \{\pi_*e_2,\pi_*e_3\}\\
 \text{are each linearly independent sets.}
\end{array}\right.
\end{equation}
(Where, as usual, $\pi\: \R^4\to \R^2$ is $\pi(x,y,p,q)=(x,y)$.)
This and the implicit function theorem imply that, if $N^2$ is any $C^1$
surface in 
$\R^4$ that is tangent to both $e_1$ and $e_2$ (or one of
the other pairs $\{e_3,e_4\}$, $\{e_1,e_4\}$ or $\{e_2,e_3\}$), then
locally $N^2$ projects over an open neighborhood of the origin in
$\R^2$. 

The number theoretic assumption\eq{a-b-ind} is what is needed to
invoke the Sternberg Normal Form (Theorem~\ref{Sternberg} in the
Appendix) and conclude there are local coordinates $\ol{x}, \ol{y},
\ol{p}, \ol{q}$ centered at the origin of $\R^4$ with
$\omega=d\ol{p}\wedge d\ol{x}+ d\ol{q}\wedge d\ol{y}$ and a function
$f(u,v)$ with $f(0,0)=0$, $f_u(0,0)=f_v(0,0)=1$ and
$H=\frac12f(\ol{p}^2-a^2\ol{x}^2,\ol{q}^2-b^2\ol{y}^2)$.  In the
coordinates $\ol{x}, \ol{y}, \ol{p}, \ol{q}$ the eigen-directions of
$L$ are given by the vectors ${\mathbf v}_{1}, {\mathbf
v}_{2},{\mathbf v}_{3},{\mathbf v}_{4}$ of\eq{vi-def} below.  So in
these coordinates we can take $e_1=\mathbf{v}_1$, $e_2=\mathbf{v}_4$,
$e_3=\mathbf{v}_3$ and $e_4=\mathbf{v}_2$.  We then can compute that
$\omega(e_1,e_2)=\omega(e_3,e_4)=\omega(e_1,e_4)=\omega(e_2,e_3)=0$,
but that $\omega(e_1,e_3), \omega(e_2,e_4)\ne 0$.  This leads to a
general existence result.

\begin{thm}\label{gen-Bruss}
Let the origin of $\R^4$ be a non-degenerate critical point of $H$ such
that the linearization $L$ of the characteristic vector field $\xi_H$
at the origin satisfies the conditions\eq{a-b-ind} and\eq{ee-proj}.
Then, with the notation above, each of the pairs $\{e_1,e_2\}$,
$\{e_3,e_4\}$, $\{e_1,e_4\}$ and $\{e_2,e_3\}$ is tangent to the jet
of a $C^\infty$ solution to $H(x,y,z_x,z_y)=0$.
\end{thm}

\begin{proof}
Consider the four submanifolds $N_{\pm,\pm}^2$ of $\R^4$ defined locally
near the origin in the coordinates $\ol{x}, \ol{y},
\ol{p}, \ol{q}$ by
$$
N_{\pm,\pm}^2:=\{ (\ol{x},\ol{y},\pm a\ol{x},\pm b\ol{y}):
(\ol{x},\ol{y})\in U\},
$$
where $U$ is a small neighborhood of $(0,0)$ in $\R^2$.  On
$N_{\pm,\pm}$ the relations $\ol{p}=\pm ax$ and $\ol{q}=\pm by$ hold.
Thus on $N_{\pm,\pm}^2$ we have
$H=\frac12f(\ol{p}^2-a^2\ol{x}^2,\ol{q}^2-b^2\ol{y}^2)=\frac12f(0,0)=0$
and therefore $N_{\pm,\pm}^2\subset \{H=0\}$.  A direct calculation
using the form of the characteristic vector field $\xi_H$ in the
coordinates $\ol{x}, \ol{y}, \ol{p}, \ol{q}$ (see\eq{xi-H} below)
shows that $\xi_H$ is tangent to $N_{\pm,\pm}^2$.  Thus
$N^2_{\pm,\pm}$ is an invariant Lagrangian surface.  
We also have
\begin{align*}
T(N^2_{+,+})_0&=\Span(e_1,e_2),\quad 
T(N^2_{-,-})_0=\Span(e_3,e_4),\\
T(N^2_{+,-})_0&=\Span(e_1,e_4),\quad
T(N^2_{-,+})_0=\Span(e_2,e_3).
\end{align*}
Therefore the assumption\eq{ee-proj} implies that locally near
$(0,0,0,0)$ each of the surfaces $N^2_{\pm,\pm}$ projects onto an open
neighborhood of $(0,0)$ in $\R^2$. Consequently
Proposition~\ref{prop:jet} implies that $N^2_{\pm,\pm}$ is the jet of
a $C^\infty$ solution to $H(x,y,z_x,z_y)=0$.
\end{proof}

If $N^2\subset \R^4$ is the jet of a solution to $H(x,y,z_x,z_y)=0$
that passes through $(0,0,0,0)$,
 then $N^2$ is an invariant Lagrangian
surface in $\R^4$.  The invariance of $N^2$ under the flow of the
characteristic vector field $\xi_H$ implies that $T(N^2)_0$ is
invariant under the linearization $L$ of $\xi_H$.  This, coupled with
the fact that $\omega$ vanishes on $T(N^2)_0$, implies that $T(N^2)_0$
is one of the four subspaces $\Span(e_1,e_2)$, $\Span(e_3,e_4)$,
$\Span(e_1,e_4)$, or $\Span(e_2,e_3)$ of $T(\R^4)_0$.  

The subspace $\Span(e_1,e_2)$ is tangent to $N^2_{+,+}$ at the origin.
This is the unstable
submanifold of $\xi_H$ at $(0,0,0,0)$ and is unique
in the sense that if $\widetilde{N}^2$ is a $C^1$ invariant
submanifold for $\xi_H$ with $T(\widetilde{N}^2)_0=\Span(e_1,e_2)$
then $\widetilde{N}^2=N_{+,+}^2$ in a neighborhood of the origin.  To
see this note that the restriction $\xi_H\big|_{\widetilde{N}^2}$ has
a source at the origin (as the eigenvalues of the restriction
$L\big|_{\widetilde{N}^2}$ are $a,b>0$) and so for all points $P$ of
$\widetilde{N}^2$ sufficiently near the origin $\lim_{t\to
-\infty}\Phi_t(P)=(0,0,0,0)$.  But this is exactly the condition that
$P$ be on the unstable submanifold.  Thus $\widetilde{N}^2=N_{+,+}^2$
as claimed.  There is a corresponding uniqueness statement for the
stable
submanifold $N^2_{-,-}$.

Call a solution $z$ to $H(x,y,z_x,z_y)=0$ with $z_x(0,0)=z_y(0,0)=0$
a {\bi solution with unstable jet\/} iff its jet near the origin is
the unstable submanifold $N^2_{+,+}$ for $\xi_H$.  Likewise, $z$ is a
{\bi solution with stable jet\/} iff its jet near the origin is the
stable submanifold $N^2_{-,-}$ for $\xi_H$.   Our discussion leads at
once to a uniqueness result.

\begin{thm}\label{stable-unique}
If the hypotheses of Theorem~\ref{gen-Bruss} hold and if $z$ is a solution to
$H(x,y,z_x,z_y)=0$ with $z(0,0)=z_x(0,0)=z_y(0,0)=0$ of either
stable or unstable type, then $z$ is unique in the sense that if
$\tilde{z}$ is any $C^2$ solution with
$\tilde{z}(0,0)=\tilde{z}_x(0,0)=\tilde{z}_y(0,0)=0$, such that the
jet of $\tilde{z}$ is tangent to the jet of $z$ at $(0,0,0,0)$, then
$\tilde{z}=z$ in a neighborhood of the origin.  (The tangency
condition of the jets is equivalent to
$\tilde{z}_{xx}(0,0)=z_{xx}(0,0)$, $\tilde{z}_{xy}(0,0)=z_{xy}(0,0)$
and $\tilde{z}_{yy}(0,0)=z_{yy}(0,0)$.)\qed
\end{thm}

Call a solution to $H(x,y,z_x,z_y)=0$ with
$z(0,0)=z_x(0,0)=z_y(0,0)=0$ whose jet at the origin is not the
stable or unstable submanifold of the characteristic vector field
$\xi_H$ a {\bi saddle solution\/}. The motivation of this
terminology is that for the eikonal equation
$z_x^2+z_y^2=a^2x^2+b^2y^2$ with $a,b>0$ the solution with unstable
jet is the convex solution $z=\frac12(ax^2+by^2)$ and the solution with
stable jet is the concave solution $z=-\frac12(ax^2+by^2)$.  The other
two obvious solutions $z=\pm \frac12(ax^2-by^2)$
(which are not unique by Theorem~\ref{saddle})
have saddle points at the origin.

We now wish to investigate in detail the lack of uniqueness of saddle
solutions.  To do this it is more convenient to work in the
coordinates $\ol{x}, \ol{y}, \ol{p}, \ol{q}$.  To simplify notation we
drop the bars and just write ${x}, {y}, {p}, {q}$ but keep
in mind
that these are not the original symplectic coordinates on $\R^4$ and
that once we have constructed the jets of solutions we have to
translate these results back to statements in the original symplectic
coordinates.  With this in mind, we assume that $H$ has form
\begin{equation}\label{H:form}
H=\frac12f({p}^2-a^2{x}^2,{q}^2-b^2{y}^2),
\end{equation}
where
$$
f(0,0)=0,\quad f_{u}(0,0)=f_{v}(0,0)=1,
$$
and
$$
\omega=d{p}\wedge d{x}+ d{q}\wedge d{y}.
$$
Then the characteristic vector field is
\begin{equation}\label{xi-H}
\xi_H=f_up\frac{\f}{\f x} + f_vq\frac{\f}{\f y} + a^2f_ux\frac{\f}{\f p}
	+ b^2f_vy\frac{\f}{\f q}.
\end{equation}
The integral curves of this vector field satisfy the differential
equation
\begin{equation}\label{L-S} 
\frac{d}{dt} 
\[ \begin{matrix} x \\ y \\ p \\ q \end{matrix} \] =
\[ \begin{matrix} 0 & 0 & f_u & 0 \\ 0 & 0 & 0 & f_v \\
f_ua^{2} & 0 & 0 & 0 \\ 0 & f_vb^{2} & 0 & 0 \end{matrix} \]
\[ \begin{matrix} x \\ y \\ p \\ q \end{matrix} \].
\end{equation}
This is not a linear system as the functions $f_u$ and $f_v$ depend on
$x$, $y$, $p$ and $q$. However we will show that because of its special
structure it can be treated almost as if it were linear.  Letting $L$
be the matrix of this system 
we find that its eigenvectors are, \emph{cf.} (\ref{e-vectors}),
\begin{equation}\label{vi-def}
{{\mathbf v}}_{1}=\[\begin{matrix} 1\\ 0\\ a\\ 0\end{matrix}\],\quad
{{\mathbf v}}_{2}=\[\begin{matrix} 0\\ 1\\ 0\\ -b\end{matrix}\],\quad
{{\mathbf v}}_{3}=\[\begin{matrix} 1\\ 0\\ -a\\ 0\end{matrix}\],\quad
{{\mathbf v}}_{4}=\[\begin{matrix} 0\\ 1\\ 0\\ b \end{matrix}\],
\end{equation}
and that
$$
L{{\mathbf v}}_{1}=f_ua{{\mathbf v}}_{1},\quad L{{\mathbf v}}_2=-f_vb
{{\mathbf v}}_2,\quad L{{\mathbf v}}_3=-f_ua{{\mathbf v}}_3,
\quad L{{\mathbf v}}_4=f_vb{{\mathbf v}}_4.
$$
Thus $L$ has a basis of eigenvectors that are independent of the
variables $x$, $y$, $p$ and $q$.  For use in the statement and proof
of Theorem~\ref{main:saddle} we define the two curves
\begin{equation}\label{ai-def}
\mathbf{a}_1(t)=t\mathbf{v}_1,\quad \mathbf{a}_2(t)=t\mathbf{v}_2,
\end{equation}
which are just the first two coordinate axes for the basis
$\mathbf{v}_1, \mathbf{v}_2, \mathbf{v}_3, \mathbf{v}_4$.  What is
special about these curves is that while we will construct many
surfaces that are jets of solutions and tangent to
$\Span(\mathbf{v}_1,\mathbf{v}_2)$, these curves will lie on all of
these surfaces.

Let $w_1, w_2, w_3, w_4$ be coordinates on $\R^4$ such that
$$
\[\begin{matrix} x\\ y\\ p\\ q\end{matrix}\]
=w_1\bv_1+w_2\bv_2+ w_3\bv_3+w_4\bv_4.
$$
Then 
\begin{equation}\label{w<->p-q}
\left\{
\begin{array}{ll}
w_1=\dfrac12\(x+p/a\),&\quad w_2=\dfrac12\(y-q/b\),\\ &\\
w_3=\dfrac12\(x-p/a\),&\quad w_4=\dfrac12\(y+q/b\).
\end{array}\right.
\end{equation}
In these coordinates the differential equations for the characteristics
become
$$
\begin{array}{ll}
w_1'(t)=f_uaw_1(t),&\quad w_2'(t)=-f_vbw_2(t),\\ & \\
w_3'(t)=-f_uaw_3(t),&\quad w_4'(t)=f_vbw_4'(t).
\end{array}
$$
Recall that $f_u(0,0)=f_v(0,0)=1$.  Therefore given $\e>0$ there is a
$\delta>0$ so that if $B(\delta):=\{(w_1,w_2,w_3,w_4):
\sum_{i=1}^4w_i^2\le \delta^2\}$ is the closed ball of radius $\delta$
at the origin, then $|f_u-1|, |f_v-1|<\e\le 1/2$ in $B(\delta)$.  We
now use that if a function $w(t)$ satisfies a differential equation
$w'(t)=c\,h(t,w(t))w(t)$ on an interval $[-\kappa,\kappa]$, where $c$ is a constant
and $|h(t,w(t))-1|\le \e$, then on $[-\kappa,\kappa]$,
$w(t)=w(0)e^{c\theta(t)t}$ for a function $\gt(t)$ that satisfies
$|1-\gt|\le \e$.  The differential equations for the characteristics
are all of this form, so if $(w_1(0),w_2(0),w_3(0),w_4(0))\in
B(\delta)$ and $t$ is so that
$(w_1(\tau),w_2(\tau),w_3(\tau),w_4(\tau))\in B(\delta)$ for $\tau$
between $0$ and $t$, then
$$
\begin{array}{ll}
w_1(t)=w_1(0)e^{a\gt_1t},&\quad w_2(t)=w_2(0)e^{-b\gt_2t},\\ &\\
w_3(t)=w_3(0)e^{-a\gt_3t},&\quad w_4(t)=w_4(0)e^{b\gt_4t},
\end{array}
$$
where
$$
|1-\gt_i|\le \e \le \frac12 \quad \text{for} \quad i=1,2,3,4.
$$

Let $\phi_+,\phi_- \:\R\to \R$ be two functions so that
\begin{enumerate}
\item $\phi_\pm$ are $C^k$ functions for some $k\ge1$ and
\item $\phi_\pm$ vanishe to order $l$ at the origin for some $l\le k$.
\end{enumerate}

\begin{lemma} \label{lemma:gamma} 
There are unique functions $\psi_{\pm}$ defined in a neighborhood of
$0$ so that if $\gamma_{\pm}$ are the curves
$$
\gamma_{\pm}(s):=s\bv_1\pm s\bv_2+\phi_{\pm}(s)\bv_3+\psi_{\pm}(s)\bv_4,
$$
then
$$
H(\gamma_{\pm}(s))\equiv 0.
$$
Moreover,
\begin{enumerate}
\item $\psi_\pm$ are $C^k$ functions and
\item $\psi_\pm$ vanish to the same order $l\le k$ at the origin
that $\phi_\pm$ do.  As above, 
this means there are
$\ol{\psi}_\pm\:\R\to \R$ so that 
\begin{equation}\label{psi-bar}
\psi_\pm(s)=\ol{\psi}_\pm(s)s^l
\end{equation}
where $\ol{\psi}_\pm(s)$ are $C^k$ functions on $\R\setminus\{0\}$
and the derivatives $\frac{d^j}{ds^j}\ol{\psi}_\pm(s)$ are bounded on
$[-1,1]\setminus \{0\}$ for $0\le j\le l$.
\end{enumerate}
(In the case $H=\frac12(p^2+q^2-a^2x^2-b^2y^2)$ we have
$\psi_{\pm}(s)=\mp({a^2}/{b^2})\phi_{\pm}(s)$.)
\end{lemma}

\begin{proof}
We will show how to find $\psi_+$, the derivation for $\psi_-$ being
similar.  In terms of the coordinates $w_1,w_2,w_3,w_4$ we are looking
for $\psi_+(s)$ so that if $w_1=s$, $w_2=s$, $w_3=\phi_+(s)$, and
$w_4=\psi_+(s)$, then $\gamma_+(s)=w_1\bv_1+w_2\bv_2+w_3\bv_3+w_4\bv_4$
is the required curve.  Writing $H$ in terms of the coordinates
$x,y,p,q$ (related to the coordinates $w_1,w_2,w_3,w_4$ by
(\ref{w<->p-q})) we have $H=\frac12f(p^2-a^2x^2, q^2-b^2y^2)$. Then we
want $\psi_+$ so that
$H(\gamma_+(s))=\frac12f(-4a^2s\phi_+(s),-4b^2s\psi_+(s))=0$.
Recalling the assumptions on $f$ (\emph{i.e.} $f(0,0)=0$ and
$f_u(0,0)=f_v(0,0)=1$) we have by the implicit function theorem that
there is a smooth function $\tilde{v}(\cdot)$ defined in a
neighborhood $U$ of $0$ with $\tilde{v}(0)=0$, $\tilde{v}'(0)=-1$ so
that $f(u,\tilde{v}(u))\equiv 0$.  Then near $(0,0)$ we have
$\{(u,v):f(u,v)=0\}=\{(u,\tilde{v}(u)):u\in U\}$.  Therefore near
$s=0$ we are trying to solve $-4b^2s\psi_+(s)=v(-4a^2s\phi_+(s))$ for
$\psi_+(s)$.  As $\tilde{v}(u)$ vanishes at the origin it can be
expressed as $\tilde{v}(u)=u\ol{v}(u)$, where $\ol{v}$ is a smooth
function and $\tilde{v}'(0)=-1$ implies $\ol{v}(0)=-1$. Therefore
$\psi_+(s)=(a^2/b^2)\phi_+(s)\ol{v}(-4a^2s\phi_+(s))$.  Note that if
$\phi_+(s)$ vanishes to order~$l$ at $s=0$, then this formula makes it
clear that the same is true for $\psi_+(s)$.  Finally, if
$H=\frac12(p^2+q^2-a^2x^2-b^2y^2)$ then $f(u,v)=u+v$ and so
$\tilde{v}(u)=-u$, $\ol{v}(u)=-1$ which implies
$\psi_+(s)=-(a^2/b^2)\phi_+(s)$.
\end{proof}

Assume that $\phi_\pm$ is defined on the interval $(-s_0,s_0)$ and let
$\Phi_t$ be the flow of $\xi_H$.  Then define $F_{\pm}\:
(-s_0,s_0)\times \R\to
\R^4$ by
\begin{equation}\label{gen-F-def}
F_\pm(s,t):= \Phi_t(\gamma_{\pm}(s)).
\end{equation}
Writing this as
$$
F_\pm(s,t)=w_1(s,t)\bv_1+w_2(s,t)\bv_2+w_3(s,t)\bv_3+w_4(s,t)\bv_4
$$
from the discussion of differential equations for the characteristics
we have
\begin{equation}\label{wi-form}
\left\{
\begin{array}{ll}
w_1(s,t)&=se^{a\gt_1t},\quad \quad\quad\  w_2(s,t)=\pm s
e^{-b\gt_2t},\\ &\\
w_3(s,t)&=\phi_{\pm}(s)e^{-a\gt_3t},\quad w_4(s,t)=\psi_{\pm}(s)e^{b\gt_4t},
\end{array}\right.
\end{equation}
where $|\gt_i-1|\le \e\le 1/2$.  We are now going use $u=w_1(s,t)$ and
$v=w_2(s,t)$ as new variables and express $w_3(s,t)$ and $w_4(s,t)$ in
terms of these variables.  As $w_2(s,t)=\pm se^{-b\gt_2t}$ this
calculation breaks up corresponding to the choice of signs.  To
simplify notation we will do the calculation in the first quadrant of
the $uv$-plane, which corresponds to choosing $+$ in the definition of
$w_2$, using the functions $\phi_+$, $\psi_+$, and restricting to
$s>0$. The calculations in the other quadrants are identical.  In the
first quadrant
$$
u=w_1(s,t)=se^{a\gt_1t}, \quad v=w_2(s,t)=se^{-b\gt_2t},
$$
with $s$, $u$, and $v$ positive.  Solving for $s$ and $t$ we get
$$
s=u^{b/(a+b)}v^{a/(a+b)}e^{ab(\gt_2-\gt_1)t}, \quad
t=\frac{1}{a\gt_1+b\gt_2} \ln\(\frac{u}{v}\).
$$
Thus
$$
e^t=\(\frac{u}{v}\)^{1/(a\gt_1+b\gt_2)}=\(\frac{u}{v}\)^{(1/(a+b))+\rho_1}
$$
where
$$
\rho_1=\frac{1}{a\gt_1+b\gt_2}-\frac{1}{a+b}=
	\frac{(1-\theta_1)a+(1-\theta_2)b}{(a+b)(a\theta_1+b\theta_2)}.
$$
As we are assuming that $|1-\theta_i|\le\e\le 1/2$ this leads to
$$
|\rho_1|\le \frac{2\e}{a+b}=:C_1\e.
$$
This implies
\begin{align*}
s&=u^{\frac{b}{a+b}}v^{\frac{a}{a+b}}(e^t)^{ab(\gt_2-\gt_1)}\\ 
	&=u^{\frac{b}{a+b}}v^{\frac{a}{a+b}}
		\(\({u}{v}^{-1}\)^{\frac{1}{a+b}+\rho_1}\)^{ab(\gt_2-\gt_1)}\\ 
	&=u^{\frac{b}{a+b}}v^{\frac{a}{a+b}}\({u}{v}^{-1}\)^{\rho_2},
\end{align*}
where
$$
\rho_2:=\left(\frac{1}{a+b}+\rho_1\right)ab(\theta_2-\theta_1).
$$
As $|\theta_1-\theta_2|\le 2\e$, $|\rho_1|\le 2\e/(a+b)\le 1/(a+b)$, 
$$
|\rho_2|\le \frac{4ab\e}{a+b}=:C_2\e.
$$

Now
\begin{align*}
w_3(s,t)&=\phi_\pm(s)e^{-a\gt_3t}\\
	&=\phi_\pm\bigl(u^{\frac{b}{a+b}}v^{\frac{a}{a+b}}
		\({u}{v}^{-1}\)^{\rho_2}\bigr)\({u}{v}^{-1}\)^{\(\frac{1}{a+b}+\rho_1\)(-a\gt		_3)}\\
&=\phi_\pm\bigl(u^{\frac{b}{a+b}}v^{\frac{a}{a+b}}
		\({u}{v}^{-1}\)^{\rho_2}\bigr)
		\({u}{v}^{-1}\)^{\(\frac{-a}{a+b}+\rho_3\)}
\end{align*}
where
$$
\rho_3:=\frac{(1-\theta_3)a}{a+b} -\rho_1a\theta_3.
$$
Using the estimates $|\rho_1|\le 2\e/(a+b)$ and $|1-\gt_3|\le\e\le 1/2$
(so that $1/2\le \gt_3 \le 3/2$) this gives
$$
|\rho_3|\le \frac{\e a}{a+b} + \frac{2\e}{a+b}\frac32a=\frac{4a\e}{a+b}=:C_3\e.
$$
Using the form of $\phi_+$ given by equation~(\ref{phi-bar}) we have
\begin{align}
w_3&=\ol{\phi}_+\bigl(u^{\frac{b}{a+b}}v^{\frac{a}{a+b}}
		\({u}{v}^{-1}\)^{\rho_2}\bigr)	
		\(u^{\frac{b}{a+b}}v^{\frac{a}{a+b}} 
		\({u}{v}^{-1}\)^{\rho_2}\)^l
		\({u}{v}^{-1}\)^{\(\frac{-a}{a+b}+\rho_3\)}\nn\\
&=\ol{\phi}_+\bigl(u^{\frac{b}{a+b}}v^{\frac{a}{a+b}}
		\({u}{v}^{-1}\)^{\rho_2}\bigr)
	u^{\frac{bl-a}{a+b}}v^{\frac{a(l+1)}{a+b}}(uv^{-1})^{l\rho_2+\rho_3}.
\label{w3-uv}
\end{align}

Likewise
\begin{align*}
w_4(s,t)&=\psi_+(s)e^{b\gt_4t}\\
	&=\psi_+\bigl(u^{\frac{b}{a+b}}v^{\frac{a}{a+b}}
		\({u}{v}^{-1}\)^{\rho_2}\bigr)
	\({u}{v}^{-1}\)^{\(\frac{1}{a+b}+\rho_1\)b\gt_4} \\
	&=\psi_+\bigl(u^{\frac{b}{a+b}}v^{\frac{a}{a+b}}
		\({u}{v}^{-1}\)^{\rho_2}\bigr)
	\({u}{v}^{-1}\)^{\(\frac{b}{a+b}+\rho_4\)},
\end{align*}
where
$$
\rho_4=\frac{(\gt_4-1)b}{a+b}+\rho_1b\gt_4.
$$
Thus
$$
|\rho_4|\le \frac{\e b}{a+b}+\frac{2\e}{a+b}\frac32b=\frac{4b\e}{a+b}=:C_4\e.
$$
Using the form of $\psi_+$ given by equation~(\ref{psi-bar}) we have
\begin{align}
w_4&=\ol{\psi}_+\bigl(u^{\frac{b}{a+b}}v^{\frac{a}{a+b}}
		\({u}{v}^{-1}\)^{\rho_2}\bigr)	
		\(u^{\frac{b}{a+b}}v^{\frac{a}{a+b}} 
		\({u}{v}^{-1}\)^{\rho_2}\)^l
		\({u}{v}^{-1}\)^{\(\frac{b}{a+b}+\rho_4\)}\nn\\
&=\ol{\psi}_+\bigl(u^{\frac{b}{a+b}}v^{\frac{a}{a+b}}
		\({u}{v}^{-1}\)^{\rho_2}\bigr)
	u^{\frac{b(l+1)}{a+b}}v^{\frac{al-b}{a+b}}(uv^{-1})^{l\rho_2+\rho_4}.
\label{w4-uv}
\end{align}

\begin{definition}\label{zero-axis}
Let $k$ be a positive integer, $m$ a positive real number, and 
$U$ a connected neighborhood of the origin in $\R^2$.  Then we
denote the class of functions defined on $U$ that are $C^k$ off of the
coordinate axes and which vanish to order $m$ along the axes by
$\caxis^{(k,m)}=\caxis^{(k,m)}(U)$.  Explicitly $h\in
\caxis^{(k,m)}(U)$ iff $h$ is $C^k$ on $U\setminus \{uv=0\}$ and for
any compact subset $V$ of $U$ there is a constant $C=C_V$ so that
$$
|h(u,v)|\le C|u|^m\quad\mbox{and}\quad |h(u,v)|\le C|v|^m
$$
hold for all $(u,v)\in V$.\qed
\end{definition}

Using this definition we now summarize and extend the calculations above.

\begin{prop}\label{prop:para}
Let $H\:\R^4\to \R$ of class $C^\infty$ be given by~(\ref{H:form}),
let $\phi_\pm\:\R\to \R$ be $C^k$ (with $1\le k\le \infty$) functions
that vanish to order $l$ at the origin (as in~(\ref{phi-bar})), and
let $\psi_\pm\:\R\to \R$ and the curves $\gamma_\pm \:\R\to \R^4$ be
given by Lemma~\ref{lemma:gamma}.  Let $w_1(s,t)\cd w_4(s,t)$ be the
solutions of the characteristic equations given by~(\ref{wi-form}) and
define $u(s,t)=w_1(s,t)$, $v(s,t)=w_2(s,t)$.  As above, we can then
solve for $s$ and $t$ in terms of $u$ and $v$.  Then define functions
on a neighborhood of the origin in $\R^2$ by
$$
g(u,v):=w_3(s,t),\quad h(u,v)= w_4(s,t),\quad\mbox{for}\ u,v \neq 0,
$$
and extend this to the coordinate axes by
$$
g(u,0)=g(0,v)=h(u,0)=h(0,v)=0.
$$
Let $m$ be any real number with 
$$
m  < \min\left\{\frac{bl-a}{a+b},\frac{al-b}{a+b}\right\}.
$$
Then there is a connected neighborhood $U$ of the origin such that
$g,h\in\caxis^{(k,m)}(U)$.  If $\phi_\pm(s)=0$ in a small neighborhood of
$s=0$, then $g$ and $h$ vanish in a neighborhood of the coordinate
axes, and are $C^k$.  The surface
$$
N^2=N^2(\phi_+,\phi_-):=\{u\mathbf{v}_1+v\mathbf{v}_2
	+g(u,v)\mathbf{v}_3+h(u,v)\mathbf{v}_4:(u,v)\in U\}
$$
is an invariant Lagrangian submanifold of $\R^4$ for some sufficiently
small neighborhood $U$ of the origin in $\R^2$.
\end{prop}

\begin{proof}
Most of this follows at once from the calculation above.  The flow
$\Phi_t$ of the characteristic vector field is $C^\infty$ and the
curves $\gamma_\pm$ are $C^k$.  Therefore $w_i(s,t)$ is a $C^k$
function of $(s,t)$ for $s\ne 0$.  So by the inverse function theorem
$s$ and $t$ are $C^k$ functions of $(u,v)$, and $g$ and $h$ are $C^k$
functions off of the coordinates axes $u=0$ and $v=0$.
%
Along the $u=0$ axis the
formula~(\ref{w3-uv}) for $w_3$ in terms of $u$ and $v$ implies that
$$
|g(u,v)|= |w_3|\le C|u|^{\frac{bl-a}{a+b}}(uv^{-1})^{|l\rho_2+\rho_3|}
$$
on any compact set.  We have
$$
|l\rho_2+\rho_3|\le (lC_2+C_3)\e.
$$
Since $\frac{bl-a}{a+b}>m$ we can, by choosing $\delta $ small enough, 
make $\e$ so small that $\frac{bl-a}{a+b}-(lC_2+C_3)\e>m$. But then 
(for a possibly larger value of $C$)
$$
|g(u,v)|= |w_3|\le C|u|^m.
$$
Similar considerations, using~(\ref{w3-uv}) and~(\ref{w4-uv}), show 
$$
|h(u,v)| \le C|u|^m, \quad |g(u,v)| \le C|v|^m,
\quad |h(u,v)| \le C|v|^m.
$$

While this only shows that this
estimate holds in a small neighborhood of
$(0,0)$, we can extend it along the $u=0$ axis as follows.  From the
construction the subset
$N^2:=\{u\bv_1+v\bv_2+g(u,v)\bv_3+h(u,v)\bv_4:(u,v)\in V\}$ is
invariant under the local flow of $\Phi_t$.  Likewise the $uv$-plane,
realized as $\R^2_{u,v}:=\{u\bv_1+v\bv_2: (u,v)\in \R^2\}$, is
invariant under the flow of $\Phi_t$, as is the line defined by $u=0$,
realized as $\R_v:=\{v\bv_2:v\in \R\}$.
The estimates above imply
that in some small neighborhood of the origin $N^2$ and $\R^2_{u,v}$
have contact to order~$m$ along $\R_v$.  But the flow $\Phi_t$ is
smooth and therefore preserves the order of contact.  Thus invariance
of $N^2$, $\R^2_{u,v}$ and $\R_v$ under $\Phi_t$ shows that $N^2$ and
$\R^2_{u,v}$ have order~$m$ contact along $\R_v$ at all points that
can be moved close to the origin by $\Phi_t$.  Similarly, we can
extend the estimates along the $v=0$ axis.  

Finally, from the construction we have that $N^2\subset \{H=0\}$ and
that the characteristic vector field $\xi_H$ is tangent to $N^2$.  Let
$X$ be a vector tangent to $N^2$, then $dH(X)=0$, because
$N^2\subset\{H=0\}$ and therefore $H$ is constant
on $N^2$.  From the definition of $\xi_H$ we have
$\omega(\xi_H,X)=-dH(X)=0$.  As $\xi_H$ is tangent to $N^2$ and
$T(N^2)$ is two dimensional, this implies that $\omega=0$ at points of
$N^2$ where $\xi_H\ne0$.  But $\xi_H$ has an isolated zero at the
origin and therefore by continuity $\omega\big|_{N^2}=0$ in some
neighborhood of the origin.  This shows that $N^2$ is an invariant 
Lagrangian surface 
near the origin and completes the proof.
\end{proof}

Recall that the coordinates $x,y,p,q$ we have been working with are not
the standard coordinates on $\R^4$, but rather the coordinates
$\ol{x},\ol{y},\ol{p},\ol{q}$ that reduced $H$ to the Sternberg normal
form.  We now need to translate Proposition~\ref{prop:para} back into
a result about solutions to the original equation
$H(x,y,z_x,z_y)=0$. However as the notion of being an invariant
Lagrangian surface is invariant under a symplectic change of
coordinates, this translation does not involve much calculation.

Let $N^2$ be the invariant Lagrangian surface from the conclusion of
Proposition~\ref{prop:para}.  Then from the construction of $N^2$ the
two curves $\mathbf{a}_1$ and $\mathbf{a_2}$ in $\R^4$ defined
by~\eq{ai-def} lie in $N^2$.  Thus $\mathbf{a}_1'(0)=\mathbf{v}_1$ and
$\mathbf{a}_2'(0)=\mathbf{v}_2$ is a basis for $T(N^2)_0$.  Then
(recalling that $\mathbf{v}_1=e_1$ and $\mathbf{v}_2=e_4$ at the
origin) the condition\eq{ee-proj} and the implicit function theorem
imply that by restricting the size of $N^2$ we can assume that it
projects over a neighborhood of $(0,0)$.  Again letting $\pi\:
\R^4\to \R^2$ be the natural projection, we then have that the two
curves $\pi\circ\mathbf{a}_1$ and $\pi\circ\mathbf{a}_2$ cross
transversely at the origin.  Let $\tilde{a}_1(t)=\pi(\mathbf{a}_1(t))$
and $\tilde{a}_2(t)=\pi(\mathbf{a}_2(t))$.  Then these curves only
depend on the curves $\mathbf{a}_1$ and $\mathbf{a}_2$ and thus only
depend on the Sternberg normal coordinates
$\ol{x},\ol{y},\ol{p},\ol{q}$ and the function $H$, but are independent
of the invariant Lagrangian surface $N^2$.  As $\tilde{a}_1$ and
$\tilde{a}_2$ are $C^\infty$ and cross transversely at the origin,
there are $C^\infty$ local coordinates $\tilde{u},\tilde{v}$ centered
at $(0,0)$ so that for $t$ near zero
$$
\tilde{u}(\tilde{a}_2(t))\equiv 0,\quad \tilde{v}(\tilde{a}_1(t))\equiv0.
$$
\begin{definition}\label{bar-zero}
Let $k$ be a positive integer, $m$ a positive real number, and $U$
a connected neighborhood of the origin in $\R^2$ contained in the
domain of the coordinates $\tilde{u},\tilde{v}$.  Then denote by
$C_{{\text{\rm $\tilde u$-$\tilde v$ axes}}}^{(k,m)}(U)$ the class of
functions defined exactly as in Definition~\ref{zero-axis} but with the
coordinates $\tilde{u}$ and $\tilde{v}$ replacing the coordinates $u$
and $v$.\qed
\end{definition}

With this definition we can now translate Proposition~\ref{prop:para}
into a statement about solutions to the differential equation
$H(x,y,z_x,z_y)=0$.

\begin{thm}\label{main:saddle}
Use the notation of Proposition~\ref{prop:para} and assume that
\eq{a-b-ind} and\eq{ee-proj} hold. Then given the $C^k$ (with $1\le
k\le \infty$) functions $\phi_\pm\:\R\to \R$ there is a solution
$z$ to $H(x,y,z_x,z_y)=0$ defined in a neighborhood $U$ of the origin
that has $N^2(\phi_+,\phi_-)$ as its jet near the origin.  Moreover,
for any real number $m$ with
$$
m  < \min\left\{\frac{bl-a}{a+b},\frac{al-b}{a+b}\right\},
$$
we have $z\in C_{{\text{\rm $\tilde u$-$\tilde v$ axes}}}^{(k+1,m+1)}(U)$.
\end{thm}

\begin{proof}
That there is a solution $z$ with $N^2=N^2(\phi_+,\phi_-)$ as jet follows
from Proposition~\ref{soln}.  As $N^2=\{(x,y,z_x,z_y): (x,y)\in U\}$
there is a gain in regularity of one derivative in going from the jet
$N^2$ to the function $z$.  As $N^2$ has regularity $\caxis^{(k,m)}$,
we see that $z$ has regularity $C_{{\text{\rm $\tilde u$-$\tilde
v$ axes}}}^{(k+1,m+1)}$.  This completes the proof.
\end{proof}

\begin{remark}\label{rml:not-unique}
Choosing $\phi_+$ and $\phi_-$ to vanish to infinite order at the
origin we see that under assumptions\eq{a-b-ind} and\eq{ee-proj} any
saddle type solution to $H(x,y,z_x,z_y)=0$ will 
have an infinite dimensional family of smooth deformations.\qed
\end{remark}

\section{Formal power series and non-existence of smooth solutions}
\label{sec:formal}

In this section we will use power series methods to investigate the
local existance and regularity of solutions to the eikonal equation,
\begin{equation}\label{eq:ikon}
z_x^2+z_y^2=h,
\end{equation}
near a point where $h=0$. We assume that $h$ is smooth and that it has
a non-degenerate Hessian at the zero we are considering. We make, if
necessary, an
affine change of coordinates of the type a translation followed by a
rotation, to bring the series expanion at this zero into the form
\begin{equation}\label{h:form}
h(x,y)=a^2x^2+b^2y^2+\sum_{m+n\ge 3}h_{m,n}x^my^n.  
\end{equation}
for some $a,b>0$.

\begin{prop}\label{lemma:sec}
If $z$ is a $C^{3}$ solution to (\ref{eq:ikon}) such that
$z(0,0) = 0$, and if $a \neq b$, then
\begin{equation}
\label{form-z}
z(x,y)=\frac12(\pm ax^2\pm by^2)+O((|x|+|y|)^3).
\end{equation}
If $a=b$ then (\ref{form-z}) still holds
after making a change of coordinates
\begin{equation}
\label{change_var}
\left[\begin{array}{c} x \\ y \end{array}\right] \mapsto
\left[\begin{array}{cc} 
\cos \xi & \sin \xi \\ -\sin \xi & \cos \xi 
\end{array}\right]
\left[\begin{array}{c} x \\ y \end{array}\right],
\end{equation}
for some angle of rotation $\xi$.
\end{prop}

\begin{proof} Suppose that $z$ has a series expansion
$z(x,y)=\alpha x+\beta y + \tfrac{1}{2} (A x^{2} + 2B xy + C y^{2}) + 
O((|x|+|y|)^{3})$. Then 
$z_x^2+z_y^2=\alpha^2+\beta^2+O(|x|+|y|)$.  But $h$ vanishes at
$(0,0)$ which, along with~(\ref{eq:ikon}), implies that $\alpha=\beta=0$.
Evaluating the second order terms of $z_{x}^{2} + z_{y}^{2}$
and equating them with the second order terms of (\ref{h:form}),
we get the equations
\begin{equation}\label{eq:coef}
A^2+B^2=a^2,\qquad B(A+C)=0,\qquad B^2+C^2=b^2.  
\end{equation}

Suppose first that $a \neq b$, then we must have $B=0$ since
otherwise the middle equation would imply that $A = -C$ which,
together with the first and the last equations, contradict the
assumption that $a \neq b$. It then follows immediately that
$A = \pm a$ and $C = \pm b$, and so (\ref{form-z}) holds. 
This holds also if $a=b$ and $B=0$, {\it i.e.} $A=\pm a$ and $C=\pm a$
with any combination of signs.
 
If $a=b$ and $B \ne 0$, then any solution of (\ref{eq:coef}) is of the
form
$$
A_{\theta} = \pm a \sqrt{1 - \theta^{2}},\quad 
B_{\theta} = a\theta,\quad 
C_{\theta} = \mp a \sqrt{1 - \theta^{2}},\quad
$$
for some $\theta \in [-1,\,1], \theta \ne 0$.
For a given $\theta$, the function 
$z = \tfrac{1}{2}(A_{\theta} x^{2} + 2B_{\theta}xy + C_{\theta} y^{2}) + \cdots$
is brought to the form (\ref{form-z}) by the change of variables
(\ref{change_var}), if $\xi$ is chosen such that
$\tan \xi = \theta/(1 + \sqrt{1 - \theta^{2}})$.
\end{proof}

As the convex and concave solutions are well understood (these are the ones
corresponding to leading terms $\pm (ax^2+by^2)/2$),
we will focus on the saddle solutions.
Thus
assume that $z$ has Taylor expansion
\begin{equation}\label{u:sad}
z(x,y)=\frac12(ax^2-by^2)+\sum_{m+n\ge 3}z_{m,n}x^my^n.
\end{equation}
(Note that, if $a=b$, we work with the rotated coordinates.) Then 
\begin{equation*}
z_x=ax+\sum_{m+n\ge3}mz_{m,n}x^{m-1}y^n  
\end{equation*}
so that
$$
z_x^2=a^2x^2+\sum_{m+n\ge 3}2maz_{m,n}x^my^n 
	+\left(\sum_{m+n\ge 3}mz_{m,n}x^{m-1}y^n\right)^2.
$$
In $\(\sum_{m+n\ge 3}mz_{m,n}x^{m-1}y^n\)^2$ we get terms  
$$
mz_{m,n}x^{m-1}y^nkz_{k,l}x^{k-1}y^l=mkz_{m,n}z_{k,l}x^{m+k-2}y^{n+l}=
mkz_{m,n}z_{k,l}x^{p}y^{q},
$$
where $p+q=m+n+k+l-2\ge k+l+1$ as $m+n\ge3$ and
likewise $p+q\ge m+n+1$ as $k+l\ge3$.  Thus there are polynomials
$P_{m,n}^{[i+j<m+n]}(z_{i,j})$ in the variables $z_{i,j}$ with $i+j<m+n$ 
such that
$$
\(\sum_{m+n\ge 3}mz_{m,n}x^{m-1}y^n\)^2=\sum_{m+n\ge
4}P_{m,n}^{[i+j<m+n]}(z_{i,j})x^my^n.
$$
Putting this all together (and setting $P_{m+n}=0$ when $m+n=3$) we get
$$
z_x^2=a^2x^2+\sum_{m+n\ge 3}(2maz_{m,n}+P_{m,n}^{[i+j<m+n]}(z_{i,j}))x^my^n.
$$

Likewise,
$$
z_y=-bx+\sum_{m+n\ge 4}nz_{m,n}x^my^{n-1},
$$
so a similar calculation gives 
$$
z_y^2=b^2y^2+\sum_{m+n\ge3}(-2nbz_{m,n}+Q_{m,n}^{[i+j<m+n]}(z_{i,j}))x^my^n.
$$
Plugging these equations into~(\ref{eq:ikon}) we get the recursion
\begin{equation}\label{eq:recur}
2(ma-nb)z_{m,n}=h_{m,n}-\left(P_{m,n}^{[i+j<m+n]}(z_{i,j})
	+Q_{m,n}^{[i+j<m+n]}(z_{i,j})\right).  
\end{equation}

\begin{thm}\label{thm:exist}
Let $h$ be a smooth function so that the Taylor expansion of $h(x,y)$ near
$(0,0)$ is of the form~(\ref{h:form}), where $a$ and $b$ are positive and
linearly independent over the rational numbers.  Then at the level of
formal power series there is a unique saddle point solution $z$
to~(\ref{eq:ikon}) with leading terms of the form~(\ref{u:sad})
\end{thm}

\begin{proof}  If $a$ and $b$ are linearly independent over the rationals,
then $ma-mb\ne0$ for all $m$ and $n$.  Thus the equation~(\ref{eq:recur})
uniquely determines all the coefficients of $z$ as required.~\end{proof}

\begin{cor}\label{cor:smooth}
With the same hypothesis as in the last theorem, there is a smooth function
$\tilde{z}$ so that $\tilde{z}_x^2+\tilde{z}_y^2-h$ vanishes to infinite
order at the origin $(0,0)$.
\end{cor}

\begin{proof}  This follows from the well known theorem of Borel
(\emph{cf.}~\cite{Borel:taylor} and \cite[Thm. 1.2.6. p.~16]{Hormander:vol1})
that given any formal power series there is a smooth (\emph{i.e}
$C^\infty$) function having the given formal power series as its
Taylor series.  Thus if $z$ is the formal power series solution given
by the theorem, then let $\tilde{z}$ be a smooth function that has~$z$
as its Taylor expansion.~\end{proof}

\begin{thm}\label{thm:non-exist}
 Let $a$ and $b$ be positive real numbers that are linearly dependent
over the rationals.  Then there is a smooth function $h$ of the
form~(\ref{h:form}) such that there is no solution $z$ of the saddle
point form~(\ref{u:sad}) to the equation~(\ref{eq:ikon}).  In
particular, if $m$ and $n$ are positive integers with $ma-nb=0$ and
$m+n\ge4$ and $c\ne0$, then for $h=a^2x^2+b^2y^2+cx^my^n$ the
equation~(\ref{eq:ikon}) will not have any solution in formal power
series (or any $C^k$ solution for $k\ge m+n$.)
\end{thm}

\begin{proof} If $a$ and $b$ are positive and dependent over the rationals,
then there are positive integers $m$ and $n$ with $m+n\ge4$ so that
$ma-nb=0$.  Therefore, if there is a solution $z$ to~(\ref{eq:ikon}) of
the form (\ref{u:sad}), the equation~(\ref{eq:recur}) implies
$$
h_{m,n}=\left(P_{m,n}^{[i+j<m+n]}(z_{i,j})+Q_{m,n}^{[i+j<m+n]}(z_{i,j})\right).
$$
But then, by adding $cx^my^m$ to $h$, $c\ne0$, we get a new $h$ of the
form~(\ref{h:form})
so that there is no solution to~(\ref{eq:ikon}) in formal power
series.  For a specific example, let $h_0=a^2x^2+b^2y^2$ and
$z_0=\frac12(ax^2-by^2)$.  Then $(z_0)_x^2+(z_0)_y^2=h_0$ 
and the equation $z_x^2+z_y^2=h_0+cx^my^n$ will have
no solution in formal power series, if
$ma-nb=0$ and $c\ne0$. ~\end{proof}

\setcounter{section}{0}
\renewcommand{\thesection}{\Alph{section}}

\section{Appendix: The stable submanifold theorem and the Sternberg
normal form}

\subsection{The stable submanifold theorem.}\label{app:stable}

The following is a statement of the Stable Submanifold Theorem adapted
to our applications.  If $k\ge 1$ and $X$ is a $C^k$ vector field on a
smooth $n$-dimensional manifold $M$, we denote by $\Phi_t$ the flow of
$X$.  That is, $\Phi_t$ is the locally defined one parameter group of
diffeomorphisms of $M$ that satisfy $\frac{d}{dt}\Phi_t(P)=X(P)$.  As
$X$ is $C^k$, the function $(P,t)\mapsto \Phi_t(P)$ is also $C^k$
({\it cf.} \cite[p.~230]{Arnold:ODE} or \cite[Thm~5 p.~86]{Lang:manifolds}).
A point $P_0\in M$ is a {\bi hyperbolic critical point\/} of $X$ iff
$X(P_0)=0$ and the linearization of $X$ at $P_{0}$,
$L\: T(M)_{P_0}\to T(M)_{P_0}$, 
has no eigenvalue with zero real part.
This implies that $P_0$ is an isolated
critical point of $X$.  Letting, as usual,
$e^{tL}=\sum_{k=0}^\infty(tL)/k!$, the {\bi stable subspace\/} of $X$
at $P_0$ is the linear subspace of $T(M)_{P_0}$ defined by
$$
T_{+}(M)_{P_0}:=\{v\in T(M)_{P_0} : \lim_{t\to \infty}e^{tL}v=0\}.
$$ 
The {\bi unstable subspace\/} is likewise defined by
$$
T_{-}(M)_{P_0}:=\{v\in T(M)_{P_0} : \lim_{t\to -\infty}e^{tL}v=0\}.
$$
If no eigenvalue of $L$ has zero real part, then basic linear algebra
implies the direct sum decomposition
\begin{equation}\label{sum-subs}
T(M)_{P_0}=T_{-}(M)_{P_0}\oplus T_{+}(M)_{P_0}.
\end{equation}
When all the eigenvalues of $L$ are real and distinct (which is the
case for most of the applications in this paper) then $T_+(M)_{P_0}$
is the span of the eigenvectors corresponding to the negative
eigenvalues and $T_-(M)_{P_0}$ is the span of the eigenvectors
corresponding to the positive eigenvalues.

\begin{thm}[Stable Submanifold Theorem]\label{stable-man}
If $k\ge 1$ and $P_0$ is a hyperbolic critical point of the $C^k$
vector field $X$, then there is a connected open neighborhood $U$ of
$P_0$ such that
$$
N_+:=\{ P\in U : \lim_{t\to +\infty}\Phi_t(P)=P_0\}
$$
(the {\bi stable submanifold\/} of $X$ at $P_0$) and
$$
N_-:=\{ P\in U : \lim_{t\to -\infty}\Phi_t(P)=P_0\}
$$
(the {\bi unstable submanifold\/} of $X$ at $P_0$) are both embedded
$C^k$ submanifolds of $U$ .  The tangent space to $N_\pm$
at $P_0$ is
$T_\pm(M)_{P_0}$ so $\dim N=\dim T_\pm(M)_{P_0}$ and
therefore~(\ref{sum-subs}) implies that $N_+$ and $N_-$ intersect
transversely and that $\dim N_++\dim N_-=\dim M$.\qed
\end{thm}

This was originally proven in various degrees of generality by
Hadamard \cite{Hadamard:stable},
Liapunov~\cite{Lyapunov:stabile-trans} and
Perron~\cite{Perron:stable}.  Modern presentations of the proof can be
found in~\cite{Hirsch-Pugh-Shub} and~\cite[Chap.~6]{Shub:stable}.

\subsection{The Sternberg normal form.}\label{sec:Sternberg}

Let $H$ be a smooth function on $\R^4$ and $\xi_H$ the characteristic
vector field of $H$ as given by (\ref{char:form}).  Let $P$ be a
critical point of $H$ (and thus a rest point of $\xi_H$) and let $L$
be the linearization of $\xi_H$ at $P$.  Then the matrix of $L$ is
given by equation~(\ref{L-def}).  If the eigenvalues of $L$ are real
and nonzero, then the formulas~(\ref{eig-val}) imply that they are of
the from $a,b,-a,-b$, for some $a,b>0$.  With this  in mind, we can give
Sternberg's normal form for a Hamiltonian near a critical point.

\begin{thm}[Sternberg~\cite{Sternberg:local3}]\label{Sternberg} 
Let $H$ be a smooth function on $\R^4$ and $P$ a critical point of $H$
with $H(P)=0$.  Assume that the eigenvalues of the linearization $L$
of $\xi_H$ at $P$ are $a,b,-a,-b$ where $a,b>0$ and that
\begin{equation}\label{Q-cond}
\text{$a$ and $b$ are linearly independent over the rational numbers.}
\end{equation}
Then there are local coordinates $\ol{x}, \ol{y}, \ol{p}, \ol{q}$
centered at $P$ and a smooth function $f(u,v)$ of two variables so that
\begin{enumerate}
\item $\ol{x}, \ol{y}, \ol{p}, \ol{q}$ are symplectic coordinates.
That is, if $\omega$ is the symplectic form then
$$
\omega=d\ol{p}\wedge d\ol{x}+ d\ol{q}\wedge d\ol{y}.
$$
\item The function $f$ satisfies $f(u,v)=u+v+O(u^2+v^2)$ or more formally:
$$
f(0,0)=0,\quad f_{u}(0,0)=f_{v}(0,0)=1,
$$ 
\item In these coordinates $H$ is given by
$$
H=\frac12f(\ol{p}^2-a^2\ol{x}^2,\ol{q}^2-b^2\ol{y}^2).
$$
\end{enumerate}
\end{thm}

This follows from Thm. 9 in Sternberg~\cite[p.~603]
{Sternberg:local3} applied to the transformations given by the
flow of a Hamiltonian vector field, using the considerations from
Section~7 in Sternberg~\cite[p.~818 -~819]{Sternberg:local}. 

\vspace{0.1in}

\begin{center}
\textsc{Acknowledgments} 
\end{center}

\vspace{0.01in}

We have benefited from numerous discussions with
Bj\"orn Jawerth, and in fact this paper evolved out of Bj\"orn's
seminars on the shape from shading problem.


\providecommand{\bysame}{\leavevmode\hbox to3em{\hrulefill}\thinspace}

\end{document}